\documentclass[12pt]{article}
\usepackage{amsmath}
\usepackage[usenames]{color}
\usepackage{mathrsfs}
\usepackage{amsfonts}
\usepackage{amssymb,amsmath}
\usepackage{CJK}
\usepackage{cite}
\usepackage{cases}
\usepackage{amsthm}

\def\bs{\backslash}
\def\d{\delta}
\def\cl{\centerline}

\def\al{\alpha}
\def\vs{\vspace*}
\def\L{\mathscr{W}_{L}(\Gamma){}}
\def\Z{\mathbb{Z}}
\def\C{\mathbb{C}}
\def\F{\mathbb{C}}
\def\QED{\hfill$\Box$}

\numberwithin{equation}{section}
\newtheorem{theo}{Theorem}[section]

\newtheorem{clai}{Claim}

\def\ptl{\partial}
\def\G{\Gamma}
\begin{document}
\begin{CJK*}{GBK}{song}

\begin{center}
{\bf\large Structures of not-finitely graded Lie algebras\\ related to generalized Virasoro algebras\,$^*$}
\footnote {$^*\,$Supported by NSF grant no.~11371278, 11001200 and 11101269 of China, the Fundamental Research Funds for the Central Universities, the grant no.~12XD1405000 of Shanghai Municipal Science and Technology Commission.
$^\dag\,$Supported by China Postdoctoral Science Foundation grant no. 2013M540709.

Corresponding author: J. Han (jzzhan@mail.ustc.edu.cn).
}
\end{center}

\cl{Qiufan Chen$^{\,*}$, Jianzhi Han$^{\,\dag}$, Yucai Su$^{\,*}$}

\cl{\small $^{\,*}$Department of Mathematics, Tongji University, Shanghai
200092, China}

\cl{\small $^{\,\dag}$Deparment of Mathematics, Sichuan University, Chengdu, Sichuan
610064, China}

\vs{8pt}

{\small
\parskip .005 truein
\baselineskip 3pt \lineskip 3pt

\noindent{{\bf Abstract:}  In this paper, we study the structure theory of a class of not-finitely graded Lie algebras related to generalized Virasoro algebras. In particular,
 the derivation algebras, the
automorphism groups and the second cohomology groups
 of these Lie algebras
are determined. \vs{5pt}

}

}
\parskip .001 truein\baselineskip 6pt \lineskip 6pt
\section{Introduction}
\setcounter{equation}{0}
As it is well-known that the Virasoro
algebra is an important object in mathematics and  physics, whose
 theory has been widely used in literature (e.g., \cite{K,KR}). Various generalizations of the  Virasoro algebra have been studied by several
authors (e.g., \cite{KR,PZ,SZ1,
SZ3,
S3,SXZ,WWY,X}). In this paper, we consider the following Lie algebras $W(\G)$:
Let $\Gamma$ be any nontrivial additive subgroup of $\C$, and $\C[\Gamma\times\Z_+]$  the semigroup algebra of $\Gamma\times\Z_+$
with basis $\{x^{\al,i}:=x^\al t^i\,|\,\alpha\in \Gamma,\,i\in\Z_+\}$ and product $x^{\al,i}x^{\beta,j}=x^{\al+\beta,i+j}$. Let $\ptl_x,\,\ptl_t$ be the derivations of $\F[\G\times\Z_+]$ defined by
$\ptl_x(x^{\al,i})=\al x^{\al,i}$, $\ptl_t(x^{\al,i})=i x^{\al,i-1}$ for $\al\in\G,\,i\in\Z_+$. Denote $\ptl=\ptl_x+t^2\ptl_t$. Then the Lie algebra $W(\G)$ is $\F[\G\times\Z_+]\ptl$ with basis $\{L_{\al,i}:=x^{\al,i}\ptl\,|\,\al\in\G,\,i\in\Z_+\}$ and \vspace*{-7pt}relation
\begin{equation}\label{algebra-rela1}[L_{\al,i},L_{\beta,j}]=
(\beta-\al)L_{\al+\beta,i+j}+(j-i)L_{\al+\beta,i+j+1},\end{equation}
for $\al,\beta\in\G,\ i,j\in\Z_+.$
One can realize the Lie algebra $W(\G)$ as follows: Let ${\cal A}:=C^\infty_{[0,+\infty)}$ be the algebra consisting of smooth functions on variable $t$ in the interval $[0,+\infty)$, which becomes a Lie algebra under bracket $[f,g]=
fg'-f'g$ for $f=f(t),\,g=g(t)\in{\cal A}$, where the prime stands for the derivative $\frac d{dt}$.
Then $W(\G)$ is the Lie subalgebra of ${\cal A}$ consisting of smooth
 functions $L_{\al,i}(t)=-\frac{e^{-\al t}}{(1+t)^i}$, $a\in\G,i\in\Z_+$.
 Thus, $W(\G)$ appears very naturally. One can observe
that this algebra
looks very similar to the Lie algebra ${\cal W}=W(0,1,0;\G)$ of Witt type studied
in \cite{SXZ}, which has the same basis with  \vspace*{-7pt}bracket $$[L_{\al,i},L_{\beta,j}]=(\al-\beta)L_{\al+\beta,i+j}+(j-i)L_{\al+\beta,i+j-1}.$$
However, $W(\G)$ has the following significant different
features from that of $\cal W$:\begin{itemize}\parskip-3pt\item[(i)]The
$\cal W$ is simple, but $W(\G)$ has infinitely many ideals.
\item[(ii)] The $\cal W$ has an ad-locally finite element $L_{0,0}$, but $W(\G)$ does not have any nonzero ad-locally finite element.
\item[(iii)] The  $\cal W$ has a {\it finitely graded filtration} in the sense that there exists a filtration
$0\subset{\cal W}^{(0)}\subset{\cal W}^{(1)}\subset\cdots,$ satisfying
$[{\cal W}^{(i)},{\cal W}^{(j)}]\subseteq{\cal W}^{(i+j)}$ for all $i,j\in\Z_+$ and
each ${\cal W}^{(i)}$ is  finitely graded, i.e., there exists some abelian group $G$ which is independent of $i$ (one can simply choose $G=\G$ in this case) such that
$${\cal W}^{(i)}=
\raisebox{-5pt}{${}^{\, \, \displaystyle\oplus}_{g\in G}$}
{\cal W}^{(i)}_g,\ \ \
[{\cal W}^{(i)}_g,{\cal W}^{(j)}_h]\subseteq {\cal W}^{(i+j)}_{g+h},\ \ \
{\rm dim\,}{\cal W}^{(i)}_g<\infty,\mbox{ \ and \ }{\cal W}_0^{(0)}\neq 0,$$ for all $g,h\in G,\,i,j\in\Z_+$ (cf.~\eqref{Fi-g}). On the other hand, $W(\G)$ does not have any such finitely graded filtration for any abelian group. (See the first part of the proof of Theorem \ref{theo1}(2).
One can also see this directly from the following: Fix $0\neq x\in W_0^{(0)}$, take any $\beta\notin{\rm Supp}\, {x}\cap \G$ (cf. \eqref{GGG-x} below) and assume that $L_{\beta,1}\in \oplus_{g\in I} W_{g}^{(j)}$ for some $j\in\Z_+$ and some finite subset $I$ of $G$. Then $\{{\rm ad}_x^k L_{\beta,1}\,|\,k\in\Z_+\}$ is an independent subset of $\oplus_{g\in I}W^{(j)}_g$, which forces some space $W_g^{(j)}$ to be infinite dimensional.)
\end{itemize}

The Lie algebra  $W(\G)$ has
the unique universal 
central extension $\hat W(\G)=W(\G)\oplus\C C$ with one dimensional center $\C C$ and  \vspace*{-7pt}relation (cf.~Theorem \ref{theo4})
\begin{equation}\label{algebra-rela2}[L_{\al,i},L_{\beta,j}]=(\beta-\al)L_{\al+\beta,i+j}+(j-i)L_{\al+\beta,i+j+1}
+\d_{\al+\beta,0}\d_{i+j,0}\frac{\al^3-\al}{12}C,\end{equation}
for $\al,\beta\in\G,\, i,j\in\Z_+.$ Then $\hat W(\G)$ contains the (generalized) Virasoro
 \vspace*{-4pt}algebra \begin{equation}\label{Vir(G)}{\rm Vir}(\G)={\rm span}\{L_{\al,0},C\,|\,\al\in\G\}.\end{equation}
We refer to $\hat W(\G)$  as a
{\it not-finitely graded} ({\it generalized}) {\it Virasoro algebra} (see Theorem \ref{theo1}(1)). Not-finitely graded Lie algebras are important objects
in Lie theory,  whose structure and representation theories are subjects of studies with more challenge
 than that of finitely graded Lie algebras.
 The Lie algebra $\hat W(\G)$ is interesting to us in another aspect that
it is also closely related to Lie algebras of Block type
 (e.g., \cite{S4,S5
 }), and it contains some interesting properties as stated in Theorem \ref{theo1}. Furthermore
 it also contains many interesting (finitely) $\G$-graded subquotient
algebras $\tilde W^{m,n}$ for $n{\sc\!}\ge{\sc\!} m{\sc\!}\ge{\sc\!}0$, where
\begin{equation}\label{bb-m}\tilde W^{m,n}=W^m/W^{n+1},\ \ \ W^m={\rm span}\{L_{\al,i}\,|\,\al\in\G,\,i\ge m\}\oplus\d_{m,0}\C C.\end{equation}
 For instance, $\tilde W^{0,0}$ is the (generalized) Virasoro algebra ${\rm Vir}(\G)$
(thus ${\rm Vir}(\G)$ is both a subalgebra and a quotient algebra of $\hat W(\G)$), $\tilde W^{0,1}$ is the well-known $W$-algebra $W(2,2)$ (up to different central extensions, see, e.g., \cite{ZD}).
When one studies the representation theory, it is well-known that the category of modules of a quotient Lie algebra $\tilde W^{0,n}$
is a full subcategory of the category of $\hat W(\G)$-modules. Thus a study of
$W$-modules will lead to a study of $\tilde W^{0,n}$-modules for all $n>0$.

 In this paper we shall mainly study the structure theory (namely, derivations, automorphisms, 2-cocycles), in this case the central element $C$ does not play a crucial role, thus we shall only consider $W(\G)$ instead of $\hat W(\G)$, and
we simply denote $W:=W(\G)$. We shall study the representation theory of $\hat W(\G)$ in a sequel paper.
The Lie algebra $W$ is $\G$-graded
\begin{equation}\label{a1-0}
W=\raisebox{-5pt}{${}^{\, \, \displaystyle\oplus}_{\alpha\in \Gamma}$}{W}_{\alpha},\ \ \
 {W}_{\alpha}={\rm span} \{L_{\alpha,i}\mid i\in \Z_+\}\mbox{ \ for }\al\in\Gamma.
 \end{equation}
However, it is not finitely graded. Nevertheless, due to the fact that $\Gamma$ may not be finitely generated (as a group), and so $W$ may not be finitely generated as a Lie algebra, the classical techniques (such as those in \cite{F}) cannot be directly applied to our situation here. One must employ
some new techniques in order to tackle  problems associated with not-finitely graded and not-finitely generated Lie algebras (this is also one of our motivations to present our results here). For instance, one of our strategies used in the present paper is to embed $W$  into its completed algebra $\overline{W}$ (see \eqref{Infini}) so that the determination of derivations can be done much more efficiently than that in some classical methods in the literature (e.g., \cite{WWY}).
The main results of the present paper are summarized in Theorems \ref{theo1}, \ref{theo2}, \ref{theo3} and \ref{theo4}.

Throughout the paper, we denote by $\C,\,\C^*,\, \Z,\, \Z_+,\,\G^*$ the sets of complex numbers, nonzero complex numbers, integers, nonnegative integers, nonzero elements of $\G$ respectively.

\section{Some properties of $W(\G)$}
We first study some properties of the Lie algebra $W:=W(\G)$, which will be summarized in Theorem \ref{theo1}.
First we recall some concepts.
A Lie algebra $\cal L$ is  {\it finitely graded} if
there exists an abelian group $G$ such that ${\cal L}=\oplus_{a\in G}{\cal L}_{[a]}$ is $G$-graded satisfying
\begin{equation}\label{Fi-g}
\mbox{$[{\cal L}_{[a]},{\cal L}_{[b]}]\subset {\cal L}_{[a+b]}$ \ and \ ${\rm dim\,}{\cal L}_{[a]}<\infty$ \ for $a,b\in G$.}\end{equation}
An element $x\in {\cal L}$ is {\it ${\rm ad}$-locally finite} if for any $y\in {\cal L}$ the subspace ${\rm span}\{{\rm ad}_x^i(y)\,|\,i\in\Z_+\}$ is finite-dimensional, where ${\rm ad}_x:y\mapsto[x,y]\ (y\in \cal L)$ is the {\it adjoint operator} of $x$.

By \eqref{bb-m} (we regard  $C$ as zero), we have
the following filtration of ideals,
\begin{equation}\label{filtra}
W=W^0\supset W^1\supset W^2\supset \cdots.
\end{equation}
\begin{theo}\label{theo1}\begin{itemize}\parskip-3pt\item[\rm(1)]The Lie algebra $W$ does not contain any nonzero ${\rm ad}$-locally finite element.
\item[\rm(2)]The Lie algebra $W$ is not finitely graded.
\item[\rm(3)]The
$W^n$
with $n\in\Z_+$ are all of the nonzero ideals of $W$ $($thus $W^1$ is the unique maximal ideal of $W)$. Furthermore\vspace*{-10pt},
\begin{equation}\label{Ideal-rela}
W^n={\rm ad}^m_{W^1}(W^{n-m})=[\ \,\stackrel{m}{\overbrace{\!\!\!W^1,[W^1,...,[W^1\!\!\!\!\!\!\!}}\ \ \ \ ,W^{n-m}]...]]\mbox{ \ for any \ }0\le m\le n.\end{equation}
\end{itemize}
\end{theo}\noindent{\it Proof.~}~(1)  Suppose $0\ne x=\sum_{\al\in\G,\,i\in\Z_+}a_{\al,i}L_{\al,i}$
(finite sum) for some $a_{\al,i}\in\C$ is an ad-locally finite element.
Choose a total order ``$\prec$'' on $\G$ compatible with its group structure, and define the total order on $\G\times\Z_+$ by \begin{equation}\label{order-1}
(\al,i)\prec(\beta,j)\ \ \ \Longleftrightarrow\ \ \ \al\prec\beta\mbox{ \ or \ }\al=\beta,\,i<j,\end{equation}
for $(\al,i),(\beta,j)\in\G\times\Z_+$.
Let $(\al_0,i_0)={\rm max}\{(\al,i)\,|\,a_{\al,i}\ne0\}$. Take $y=L_{0,j}$ with $j=\d_{i_0,0}$. Then for different $k$, the element ${\rm ad}^k_x(y)$ has the different highest  \vspace*{-7pt}term $$\mbox{$\prod\limits_{p=0}^{k-1}$}\Big(j-i_0+p(i_0+1)\Big)a_{\al_0,i_0}^kL_{k\al_0,j+k(i_0+1)}\ne0,$$ i.e., ${\rm ad}^k_x(y)$ for $k=0,1,2,...,$ are linearly independent, which is a contradiction  \vspace*{-7pt}with $${\rm dim\,}{\rm span}\{{\rm ad}_x^i(y)\,|\,i\in\Z_+\}<\infty.$$

(2) Suppose  there exists some abelian group $G$ such that $W=\oplus_{g\in G}W_{[g]}$ is $G$-graded satisfying \eqref{Fi-g}. Suppose $W_{[0]}\neq 0$.  Take any nonzero $x\in W_{[0]}$ and $\beta\in \G$ such that $\beta\notin {\rm Supp}x$ (cf. \eqref{GGG-x} below). Assume that $L_{\beta,1}\in \oplus_{g\in I}W_{[g]}$ for some finite subset $I$ of $G$.  Then it follows from  the proof of (1) that $\{{\rm ad}^k_x L_{\beta,1}|k\in\Z_+\}$ is an independent subset of $\oplus_{g\in I}W_{[g]}$. In particular, $\sum_{g\in I}{\rm dim}W_{[g]}=\infty$. Since $I$ is finite, there exists some $g_0\in I$ such that dim$W_{[g_0]}=\infty$, contradicting our assumption that  $W_{[g]}$ is finite dimensional for all $g\in G$.

So in the following we assume that $W_{[0]}=0$.
Choose a total order ``$\prec$'' on $G$ compatible with its group structure.  Write $L_{0,0}=\sum_{i=1}^nx_{g_i}$ with $g_n\succ g_{n-1}\succ\cdots \succ g_1$ for some $n\geq 1$, where $0\neq x_{g_i}\in W_{[g_i]}$.  Since
\begin{equation}\label{eq-alpha-------}[L_{0,0}, L_{\alpha, 0}]=\alpha L_{\alpha,0}
\end{equation}that is, $L_{0,0}$ preserves each one dimensional space $\C L_{\alpha,0}$ ($\alpha\in\G$),  $L_{0,0}$ can not lie in the positive part $\oplus_{0\prec g\in G}W_{[g]}$ or the negative part $\oplus_{0\succ g\in G}W_{[g]}$ of $W$. So  $g_n\succ 0$ and $g_1\prec 0$.

By comparing the maximal homogenous components (with respect to the $G$-gradation of $W$) of both sides of \eqref{eq-alpha-------}
 we see that the maximal homogenous component of $L_{\alpha,0}$ lies in $W_{[g_n]}$ for all $\alpha\in\Gamma$. So there exists $\lambda_{\alpha}\in\C^*$ such that $L_{\alpha,0}=\lambda_{\alpha}L_{0,0}+\sum_{g\prec g_n}y^{\alpha}_g$, where $y^\alpha_g\in W_{[g]}$. In fact, substituting $L_{\alpha,0}=\lambda_{\alpha}L_{0,0}+\sum_{g\prec g_n}y^{\alpha}_g$ in the left hand side of \eqref{eq-alpha-------} and then proceeding the analysis above give $L_{\alpha,0}=\lambda_{\alpha}L_{0,0}+\sum_{g\prec0 }y^{\alpha}_g$.
Furthermore, comparing the minimal homogenous component shows that $L_{\alpha,0}=\lambda_\alpha L_{0,0}+y^{\alpha}_{g_1}$ for some  $y^{\alpha}_{g_1}\in W_{[g_1]}$. In this case, we have
\begin{equation*}\infty={\rm dim(span}\{L_{\alpha,0}-\lambda_{\alpha}L_{0,0}| \alpha\in\Gamma\})\leq \, {\rm dim}W_{[g_1]},
\end{equation*} contradicting the finite dimension of $W_{[g_1]}$.

(3)
Suppose $A$ is a nonzero ideal, and
$0\ne x=\sum_{\al\in\G}x_\al\in A$ (finite sum) such that each homogeneous component $x_\al=\sum_{i\in\Z_+}a_{\al,i}L_{\al,i}$ (finite sum) for some $a_{\al,i}\in\C$.
Let \begin{equation}\label{GGG-x}{\rm Supp\,}x=\{\al\in\G\,|\,x_{\al}\ne0\}\mbox{ \ (called the {\it support} of $x$)}.\end{equation}
We call the size $|{\rm Supp\,}x|$ is the {\it depth} of $x$.
We want to prove the following claim by induction on the depth of $x$.
\begin{clai}There exists some basis element $L_{\beta_0,j_0}\in A$.
\end{clai}
First suppose $|{\rm Supp\,}x|=1$, i.e., ${\rm Supp\,}x=\{\al\}$ for some $\al\in\G$.
Thus $x$ is a homogeneous  element and we can simply write $x=x_\al=\sum_{k=i}^{j}a_k L_{\al,k}$ for some $i\le j$ such that $a_{i},a_{j}\ne0$. We say $a_{i}L_{\al,i},\,a_{j}L_{\al,j}$ are respectively the {\it first, last term} of $x_\al$, and $i,j$ the {\it first, last index} of $x_\al$. We denote\begin{equation}\label{length---}\ell(x)=j-i+1\mbox{ (called  the {\it length} of the homogeneous element $x$)}.\end{equation} If $\ell(x)=1$, we have the claim. Thus suppose
$\ell(x)>1$.  \vspace*{-7pt}Take \begin{equation}\label{y1first}
y=[L_{\al,j},x]=\mbox{$\sum\limits_{k=i}^{j-1}$}a_{k}(k-j)L_{2\al,j+k+1}\in A.\end{equation}
 We see that $y\ne0$ is a homogeneous element such that $\ell(y)<\ell(x)$. The claim is obtained by induction on the length $\ell(x)$.

Now suppose $|{\rm Supp\,}x|>1$. Take a nonzero homogeneous component $x_\al\ne0$ and set
$$\mbox{$y:=[x_{\al},x]=\sum\limits_{\beta\in{\rm Supp\,}x,\,\beta\ne\al}[x_\al,x_\beta]\in A$.}$$
One can easily see from \eqref{algebra-rela1} that the first term in $[x_{\al},x_\beta]$ is nonzero if $\beta\ne\al$, thus $y\ne0$ and ${\rm Supp\,}y=\{\al+\beta\,|\,\beta\in({\rm Supp\,}x)\bs\{\al\}\}$ has depth $|{\rm Supp\,}x|-1$. The claim is proved by induction on the depth of $x$.\vskip5pt

Now suppose $j_0$ is smallest such that Claim 1 holds for some $\beta_0\in\G$.
Fix any $\gamma\in\G^*$.
For any $\beta\in\G$,
we define the \vspace*{-7pt}operator
$$\theta_\beta:={\rm ad}_{L_{\beta-\gamma,0}}{\rm ad}_{L_{\gamma,0}}-2\,{\rm ad}_{L_{\beta,0}}{\rm ad}_{L_{0,0}}+
{\rm ad}_{L_{\beta+\gamma,0}}{\rm ad}_{L_{-\gamma,0}}.$$
Since $A$ is an ideal, applying $\theta_\beta$
to
$L_{\beta_0,j_0}$, we see that $-4\gamma^2L_{\beta+\beta_0,j_0}=\theta_\beta(L_{\beta_0,j_0})\in A$. From this one can deduce
that $L_{\beta,j}\in A$ for all $(\beta,j)\in\G\times\Z_+$ with $j\ge j_0$. Hence, $A\supset W^{j_0}$. If $j_0=0$, then $A=W^0=W$. So we assume $j_0\geq 1$ in the following.

Suppose $A\nsubseteq W^{j_0}$. Consider the nonzero quotient space $A/W^{j_0}$ and denote by $\bar x$ the image of $x$ in $A/W^{j_0}$ for any $x\in A$. Then we can take $0\neq \bar x=\sum_{\al\in I,\ 0\leq i\leq j_0-1}c_{\al, i}\bar L_{\al, i}$, where $I$ is a finite subset of $\G$. Set $i_0={\rm min}\{i\ |\ c_{\al,i}\neq 0, \al\in I\}$. Fix any $\beta\in \G-I$, then one can see that the image  of $y_\beta:=[L_{\beta,j_0-i_0-1}, x]$ is nonzero in $W^{j_0-1}/W^{j_0}\bigcap A/W^{j_0}$. Observe that $\bar L_{0,0}$ acts semisimply on $W^{j_0-1}/W^{j_0}$.  Repeatedly applying the action of $\bar L_{0,0}$ on $\bar y_\beta$ will yield that $L_{\al+\beta, j_0-1}\in A$ for $\al\in I$ such that $c_{\al, i_0}\neq 0$, a contradiction with the minimality of $j_0$.\QED

\def\L{W}
\section{ Derivation algebra}
Recall that a linear map $D:\L\rightarrow \L$ is a {\it derivation} of $\L$
if $D\big([x,y]\big)=[D(x),y]+[x, D(y)]$ for any $x,y\in\L$.
For any  $z\in\L$,
 the adjoint operator ${\rm ad}_z\!:\L\rightarrow\L$ is a derivation, called an {\it inner
 derivation}. Denote by
 ${\rm Der\,}{\L}$
 and ${\rm ad\,}{\L}$ the vector spaces of all derivations and inner derivations respectively.
 Then the first cohomology group
${H}^{1}(\L,\ \L)\cong{\rm Der\,}{\L}/{\rm
ad\,}{\L}$.

Let ${\rm Hom}_\Z(\Gamma,\C)$ denote the space of
group homomorphisms from $\Gamma$ to (the additive group) $\F$ (for each $\phi\in {\rm Hom}_\Z(\Gamma,\C)$, the scalar multiplication $\phi$ by $c\in\C$ is defined by $(c\phi)(\gamma)=c\phi(\gamma)$, thus ${\rm Hom}_\Z(\Gamma,\C)$ is a vector space). For each $\phi\in{\rm Hom}_\Z(\Gamma,\C)$, we can define a derivation
$D_{\phi}$ as follows,
\begin{equation}\label{d5}
D_{\phi}( L_{\al,i})=\phi(\alpha)L_{\al,i}
\ \ \ \mbox{for}\ \ \al\in\Gamma,\ i\in\Z_+.
\end{equation}
We still use ${\rm Hom}_\Z(\Gamma,\F)$ to  denote the
corresponding subspace of ${\rm Der\,}\L$. In particular, since $\phi_0:\al\mapsto\al$ is in ${\rm Hom}_\Z(\Gamma,\C)$, we have the derivation
\begin{equation}\label{D0====}D_0=D_{\phi_0}:L_{\al,i}\mapsto \al L_{\al,i}\mbox{ \ for \ }\al\in\G,\ i\in\Z_+.\end{equation}
\begin{theo}\label{theo2}
We have ${\rm Der\,}\L={\rm ad\,}\L\oplus{\rm Hom}_\Z(\Gamma,\F)$. In particular, if $\G=\Z$, we have ${\rm Der\,}W(\Z)={\rm ad\,}W(\Z)\oplus\C D_0$.\end{theo}
\noindent{\it Proof.~}~First by noting that every element in ${\rm Hom}_\Z(\Gamma,\F)$ acts locally-finitely on $W$, and by Theorem \ref{theo1}(2), we can easily prove that the sum in the theorem is direct. Since ${\rm Hom}_\Z(\Z,\F)=\C\phi_0,$ we have the second statement.

Now
\def\ol{\overline}take  $\ol{W}$ to be the space whose elements are finite sums of homogeneous elements $x_\al$ for $\al\in\G$ (i.e., each element of $\ol W$ has a finite support, cf.~\eqref{GGG-x}), and each $x_\al$ has the form $x_\al=\sum_{i\in\Z_+}a_{\al,i}L_{\al,i}$ (which can be an infinite sum, i.e., each homogeneous element can have infinite length, cf.~\eqref{length---}).
Then $\ol{W}$ becomes a Lie algebra with bracket defined by (cf.~\eqref{algebra-rela1})
\begin{equation}\label{Infini}\mbox{$\Big[\sum\limits_{\al,i}a_{\al,i}L_{\al,i},
\sum\limits_{\beta,j}b_{\beta,j}L_{\beta,j}\Big]=\sum\limits_{\al,\beta,i,j}a_{\al,i}b_{\beta,j}\Big((\beta-\al)L_{\al+\beta,i+j}
+(j-i)L_{\al+\beta,i+j+1}\Big).$}\end{equation}
Clearly $W\subset\ol W$. We call $\ol W$ the {\it completion} of $W$. Let $D\in{\rm Der\,}\L$.
Applying $D$ to $[L_{0,0}, L_{\alpha,i}]=\alpha L_{\alpha, i}+i L_{\alpha, i+1}$ one can see that $${\rm Supp}(D(L_{\alpha, i+1}))\subseteq {\rm Supp}((DL_{\alpha, i}))\cup{\rm (Supp}(D(L_{0,0})+\alpha)$$ for any $i\geq 1$. Inductively, we have $${\rm Supp}(D(L_{\alpha, i+1}))\subseteq {\rm Supp}(D(L_{\alpha, 1}))\cup{\rm (Supp}(D(L_{0,0})+\alpha)$$ for any $i\geq 1$. So $\sum_{i\in\Z_+}a_{i}D(L_{\alpha, i})\in \ol W$. Thus if we define     $$\mbox{$\ol D\Big(\sum\limits_{\al,i}a_{\al,i}L_{\al,i}\Big)=\sum\limits_{\al,i}a_{\al,i}D(L_{\al,i})$}$$  for any $\sum_{\alpha,i}a_{\alpha, i}L_{\alpha,i}\in \ol W$,  then $\ol D\in {\rm End\,}\ol\L$; Moreover, $\ol D\in {\rm Der\,}\ol\L$.
We shall prove the result in two steps: \begin{itemize}\parskip-3pt
\item[(i)]First by replacing $\ol D$ by  $\ol D-{\rm ad}_y-D_\phi$ (here we simply denote $\ol D_\phi$ by $D_\phi$) for some $y\in\ol W,\,\phi\in{\rm Hom}_\Z(\Gamma,\F)$, we want to prove $\ol D|_W=0$.\item[(ii)] Then we want to prove that in fact $y\in W$.
\end{itemize}

Suppose $\ol D(L_{0,0})\!=\!\sum_{\al,i}a_{\al,i}L_{\al,i}\in\L$ for some $a_{\al,i}\in\C$.
For any $\al\in\G$, we define $b_{\al,j}\in\C$ inductively on $j\ge0$ by  regarding $b_{\al,-1}$ as zero and
\begin{equation}\label{definedbbbb}
b_{\al,j}=\left\{\begin{array}{llll}\frac1\al(-a_{\al,j}-b_{\al,j-1}(j-1))&\mbox{if \ }\al\ne0,\\[4pt]-\frac1ja_{0,j+1}&\mbox{if \ }\al=0,\,j>0,\\[4pt]
0&\mbox{if \ }\al=j=0.
\end{array}\right.\end{equation}
Take $y=\sum_{\al,j}b_{\al,j}L_{\al,j}$. Note that $y\in\ol{W}$ (since the depth of $y$ is $|{\rm Supp\,}y|\le|{\rm Supp\,}x|$, cf.~\eqref{GGG-x}), and \eqref{definedbbbb}  \vspace*{-7pt}gives $$\begin{array}{ll}
 \ol D(L_{0,0})-{\rm ad}_y(L_{0,0})\!\!\!&=\sum\limits_{\al,j}a_{\al,j}L_{\al,j}-\sum\limits_{\al,j}\Big
 (-b_{\al,j}\al-b_{\al,j-1}(j-1)\Big)L_{\al,j}
 \\[12pt]&
=a_{0,0}L_{0,0}+a_{0,1}L_{0,1}. \end{array}$$This proves the following claim.
\setcounter{clai}{0}\begin{clai}\label{clai2}
By replacing $\ol D$ by $\ol D-{\rm ad}_y$ for some $y\in\ol W$, we can suppose \begin{equation}\label{L-00-D}
\ol D(L_{0,0})=a_{0}L_{0,0}+a_{1}L_{0,1}\mbox{ \ for some \ }a_0,a_1\in\C.\end{equation}
\end{clai}
For any $\al\in\G$, assume $\ol D(L_{\al,0})=\sum_{\beta,j}a_{\beta,j}^\al L_{\beta,j}$ for some $a_{\beta,j}^\al\in\C$.
Applying $\ol D$ to $[L_{0,0},L_{\al,0}]=\al L_{\al,0}$ and comparing the coefficients of $L_{\beta,k}$,
by induction on $k$ we obtain $a^\al_{\beta,k}=0$ for all $k\ge 0$ if $\beta\ne\al$ and $a^\al_{\al, k}=0$ for all $k\ge 1$, and we also obtain $a_0=a_1=0$.
Thus we can assume $\ol D(L_{\al,0})=b_{\al} L_{\al,0}$ for some $b_{\al}\in\C$ such that by \eqref{L-00-D}, \begin{equation}\label{b0000k}b_{0}=0.\end{equation}
Applying $\ol D$ to $[L_{\al,0},L_{\gamma,0}]=(\gamma-\al)L_{\al+\gamma,0}$, we obtain that the following holds for $\al\ne\gamma$,
\begin{equation}\label{al===nnn}b_{\al+\gamma}=b_{\al}+b_{\gamma}.\end{equation}
In particular, $b_{-\al}=-b_{\al}$ for all $\al\in\G$.
If $\al=\gamma$, we have $b_{2\al}=b_{(\al+\eta)+(\al-\eta)}=b_{\al+\eta}+b_{\al+(-\eta)}=b_{\al}+b_{\eta}+
b_{\al}+b_{-\eta}=2b_{\al}$ for any $\eta\in\G\bs\{\pm\al,0\}$. Thus \eqref{al===nnn} holds for all $\al,\gamma\in\G$, which shows that the map $\phi:\al\mapsto b_{\al}$ is an element in ${\rm Hom}_\Z(\Gamma,\F)$. If we replace $\ol D$ by $\ol D-D_\phi$ (cf.~\eqref{d5}) (note that this replacement does not affect \eqref{L-00-D}), we can suppose $b_{\al}=0$ for all $\al\in\G$.
This proves the following claim.
\begin{clai}\label{clai2-1}
Replacing $\ol D$ by $\ol D-D_\phi$ for some $D_{\phi}\in{\rm Hom}_\Z(\Gamma,\F)$, we can suppose
$\ol D(L_{\al,0})=0$ for all $\al\in\G$.
\end{clai}

  For any $D\in {\rm Der\,}\L$ and $\triangle\in\G$, define the homogeneous operator $\ol D_\triangle$ of degree $\triangle$ in the following way
$$\ol D_\triangle(\sum_{\al}u_\al)=\sum_{\al}\pi_{\al+\triangle} \ol Du_\al,$$ where $u_\al\in \ol\L_\al$ and  $\pi_{\al}: \ol\L\rightarrow \ol\L_{\al}$ is the natural projection.
Then $\ol D=\sum_{\triangle\in\G}\ol D_\triangle$ and $\ol D_{\triangle}\in{\rm Der\,}\ol\L$.  For any $\al\in\G$, suppose
\begin{equation}\label{eq-+1}
\ol D_\triangle L_{\al,1}=\sum_{j}c_{\al,j}^\triangle L_{\al+\triangle,j}\ (c_{\al,j}^\triangle\in\C).
\end{equation}
Applying $\ol D_\triangle$ to $[L_{\al, 0}, L_{\beta, 1}]-[L_{0,0}, L_{\al+\beta,1}]=-2\al L_{\al+\beta,1}$, one has
\begin{equation*}
(\beta+\triangle-\al)(c_{\al+\beta,j}^\triangle-c_{\beta,j}^\triangle)=(1-j)(c_{\al+\beta, j-1}^\triangle-c_{\beta,j-1}^\triangle)\ ({\rm we\ take}\ c^\triangle_{\al,-1}=0) .
\end{equation*}
It follows that  $c_{\al+\beta,j}^\triangle=c_{\beta,j}^\triangle$ for any $j\geq 0$ whenever $\beta+\triangle-\al\neq 0$ and $c_{2\beta+\triangle,j}^\triangle=c_{\beta,j}^\triangle$ for any $j\geq 1$. In particular, take $\beta=0$ and $\al+\beta=0$ respectively, then one can see that  $c_{\gamma, j}^\triangle=c_{0,j}^\triangle$ for any $\gamma\in\G$ and $j\geq 0$.
This means the coefficient $c_{\gamma, j}^\triangle$ is independent of the choice of $\gamma$. So in the following, for convenience,  we can drop off the symbol $\gamma$ and just write $c_j^\triangle$ instead of $c_{\gamma, j}^\triangle$. Note that $[L_{\al, 1}, L_{\beta,1}]+(\al-\beta)[L_{0,0}, L_{\al+\beta,1}]=(\al^2-\beta^2)L_{\al+\beta,1}$. By applying $\ol D_\triangle$ to this equation and comparing the coefficients of $L_{\al+\beta+\triangle, j}$, we obtain $\triangle c_{j}^\triangle=(3-j)c_{j-1}^\triangle$, which gives rise to $c_j^\triangle=0$ for any $\triangle \neq 0$ and $j\geq 0$ and $c_k^0=0$ whenever $k\neq 2$. Hence the equation (\ref{eq-+1}) can be simply written as $$\ol D_{\triangle} L_{\al, 1}=c\delta_{\triangle, 0}L_{\al, 2}\ (c\in \C).$$
Applying $\ol D_0$ to $[L_{\al,1}, L_{\beta,1}]=(\beta-\al)L_{\al+\beta,2}$, one can immediately see that $\ol D_0 L_{\al, 2}=cL_{\al, 3}$.  By applying $\ol D_0$ to $[L_{0,0}, L_{\al, 1}]-\al L_{\al,1}=L_{\al,2}$ we obtain $\ol D_0 L_{\al, 3}=3cL_{\al, 4}$. By induction, we can obtain that $\ol D_0 L_{\al, i}=ciL_{\al, i+1}$, that is to say $\ol D_0=c({\rm ad}_{L_{0,0}}-D_{0})$. Replacing $y$ by $y-cL_{0,0}$ and $\phi$ by $\phi^{'}:=\phi +c\phi_{0}$ respectively (which does not affect $\ol D(L_{\alpha,0})=0$), we get $\ol D L_{\al, 1}=0$. So far, we have proved that $\ol D L_{\al, 0}=\ol D L_{\al, 1}=0$ for all $\al\in\G$. Since $W$ is generated by $\{L_{\al,0},L_{\al,1}\,|\,\al\in\G\}$, we obtain $\ol D|_W=0$, and the first step is completed.

 The first step in fact shows that $D=\ol D|_W=({\rm ad}_y+D_\phi)|_W$.
If $y\notin W$, then some homogeneous component, say, $y_\al=\sum_{i\ge0}a_iL_{\al,i}$ (for some $a_i\in\C$), of $y$ is not in $W$. Choose $\beta\in\G^*$, and set $u=D(L_{\beta,0}),\,w=D(L_{2\beta,0})\in W$. Then the homogeneous components
$$\begin{array}{lll}u_{\al+\beta}&\!\!\!=\!\!\!&
-\mbox{$\sum\limits_{i\ge0}(a_i(\al-\beta)+(i-1)a_{i-1})L_{\al+\beta,i}+\d_{\al,0}D_{\phi^\prime} L_{\beta,0}$ (we take $a_{-1}=0$)},\\[12pt]w_{\al+2\beta}&\!\!\!=\!\!\!&
-\mbox{$\sum\limits_{i\ge0}(a_i(\al-2\beta)+(i-1)a_{i-1})L_{\al+2\beta,i}+\d_{\al,0}D_{\phi^\prime}L_{2\beta,0}$},\end{array}$$
are in $W$. Thus there exists some $N>0$ such that $a_i(\al-\beta)+(i-1)a_{i-1}=a_i(\al-2\beta)+(i-1)a_{i-1}=0$ for $i>N$. In particular $a_i=0$ for $i>N$, a contradiction with the fact that $y_\al\notin W$. This completes the proof of the second step, thus the theorem.\QED

\section{Automorphism group
}
\setcounter{theo}{0}
 Denote by ${\rm
Aut\,}{\L}$
the automorphism group of $\L$. Let $\chi(\G)$ be the set of characters of $\G$, i.e., the set of group homomorphisms $\tau:\G\to\C^*$. Set $\G^{\C^*}=\{c\in\C^*\,|\,c\G=\G\}$. We define a group structure on $\chi(\G)\times\G^{\C^*}$ by
\begin{equation}\label{group-sss}
(\tau_1,c_1)\cdot(\tau_2,c_2)=(\tau,c_1c_2),\mbox{ \ where \ }\tau:\al\mapsto\tau_1(c_2\al)\tau_2(\al)\mbox{ for }\al\in\G.
\end{equation}It turns out that the group $\chi(\G)\times\G^{\C^*}$ is just the semidirect product $\chi(\G)\rtimes\G^{\C^*}$ under the action given by $(c\tau)(\alpha)=\tau(c\alpha)$ for all $c\in\Gamma^*$, $\tau\in\chi(\G),$  $\alpha\in\Gamma$.
 We define a group homomorphism $\phi:(\tau,c)\mapsto\phi_{\tau,c}$ from $\chi(\G)\times\G^{\C^*}$ to ${\rm Aut\,}W$ such that
 $\phi_{\tau,c}$ is the automorphism of $W$ defined by\begin{equation}\label{Aususus}
 \mbox{$\phi_{\tau,c}:L_{\al,i}\mapsto\tau(\al)c^{-i-1}L_{c\al,i}$ \ for \ $\al\in\G,\,i\in\Z_+$.}\end{equation}
One can easily verify that $\phi_{\tau,c}$ is indeed an automorphism of $W$.
%
\begin{theo}\label{theo3}
We have $\phi:{\rm Aut\,}W\cong\chi(\G)\rtimes\G^{\C^*}$.
\end{theo}  \noindent{\it Proof.~}~Let $\sigma\in{\rm Aut\,}W$. Since $W^1$ is the unique maximal ideal of $W$, we have $\sigma(W^1)=W^1$. Then by \eqref{Ideal-rela}, we obtain $\sigma(W^n)=W^n$ for all $n\in\Z_+$.
In particular, $\sigma$ induces an automorphism $\bar\sigma$ of $\tilde W^{0,0}\cong{\rm Vir}(\G)$ (cf.~\eqref{Vir(G)} and \eqref{bb-m}). Therefore
$\sigma(L_{\al,0})\equiv \tau(\al) c^{-1} L_{c\al,0}\,({\rm mod\,}W^1)$ for $\al\in\G$ and some $\tau\in\chi(\G)$, $c\in\G^{\C^*}$ (e.g., \cite{WWY}), this result can also be proved directly by noting that $\C^*L_{0,0}$ is the set of nonzero ad-locally finite elements in ${\rm Vir}(\G)$). By replacing $\sigma$ by $\sigma\phi_{\tau,c}^{-1}$, we can suppose $\tau=1,\,c=1$. Thus
\begin{equation}\label{sigma-111}
\sigma(L_{\al,0})-L_{\al,0}\in W^1.\end{equation}
\setcounter{clai}{0}\begin{clai}\label{Auto-lamm}
We have $\sigma(L_{0,0})=L_{0,0}.$
\end{clai}

First suppose that there exists some $\al\in\G^*$ such that the homogeneous component $\sigma(L_{0,0})_\al$ of $\sigma(L_{0,0})$ is nonzero (we always use the same symbol with subscript $\al$ to denote its  \linebreak homogeneous
component of degree $\al$). Recall from \eqref{GGG-x} that the support of $\sigma(L_{0,0})$ is \linebreak ${\rm Supp\,}\sigma(L_{0,0})=\{\al\in\G\,|\,\sigma(L_{0,0})_\al\ne0\}$. Choose a total order of $\G$ compatible with its group structure such that
the maximal element $\al_0={\rm max\,Supp\,}\sigma(L_{0,0})$ of the support of $\sigma(L_{0,0})$ is $\succ0$. Set
\begin{equation}\label{maxi-elem-in}
\gamma_0={\rm max\,}{\rm Supp\,}\sigma(L_{2\al_0,0}).
\end{equation}
Since $2\al_0\in{\rm Supp}\,\sigma(L_{2\al_0,0})$ by  \eqref{sigma-111}, we have
 $
 \gamma_0\succeq2\al_0\succ\al_0\succ0.$ 
 Thus $$[\sigma(L_{2\al_0,0}),\sigma(L_{0,0})]_{\gamma_0+\al_0}\ne0,$$ by considering its first term (recall this notion in the proof of Claim 1 in the proof of Theorem \ref{theo1}). This shows that $\sigma(L_{2\al_0,0})_{\gamma_0+\al_0}\ne0$, a contradiction with the definition of $\gamma_0$ in \eqref{maxi-elem-in}. Thus  $\sigma(L_{0,0})=\sigma(L_{0,0})_0$ is a homogeneous element concentrated on degree zero. Applying $\sigma$ to $[L_{0,0},L_{\al,0}]=\al L_{\al,0}$, we see $\sigma(L_{\al,0})=\sigma(L_{\al,0})_{\al}$ is a homogeneous element concentrated on degree $\al$ for all $\al\in\G$.
If $\sigma(L_{0,0})\ne L_{0,0}$, by considering the last term of $[\sigma(L_{0,0}),\sigma(L_{\al,0})]$ with $\al\ne0$, we see $[\sigma(L_{0,0}),\sigma(L_{\al,0})]\ne\al\sigma(L_{\al,0})$, a contradiction. The claim is proved.\vskip5pt

Applying $\sigma$ to $[L_{0,0},L_{\al,0}]= \al L_{\al,0}$, by \eqref{sigma-111} and Claim 1, we obtain \begin{equation}\label{sigma-l-a-0}
\sigma(L_{\al,0})=L_{\al,0}\mbox{ \ for all \ }\al\in\G.\end{equation}
Now for any $\al\in\G$, we set \begin{equation}\label{y-alpha}\mbox{$
y^\al=\sigma(L_{\al,1})=\sum\limits_{\beta}y^\al_\beta,$  \ where \
$y^\al_{\beta}=\sum\limits_{j\ge0}c^\al_{\beta,j}L_{\beta,j}$,}\end{equation}
for some $c^\al_{\beta,j}\in\C$ with $c^\al_{\beta,0}=0$ (recall that $\sigma(W^1)=W^1$).
\begin{clai}\label{Annnclai}
We have $y^\al_\beta\!=\!0$ if $\beta\!\ne\!\al$, and $y^\al\!=\!y^\al_\al\!=\!\sum_{j\ge1}c_jL_{\al,j}$ for some $c_j\!\in\!\C$ with $c_1\!\ne\!0$.\end{clai}

Assume $y^{\al_0}$ is not concentrated on its $\al_0$-component, i.e., $y^{\al_0}\ne y^{\al_0}_{\al_0}$ for some $\al_0\in\G$. 
Let $\beta_0$ be the element in $\G^*$ (we choose an order of $\G$ such that $\beta_0\succ0$) such that
\begin{equation}\label{max-be-----}
\beta_0\succeq{\rm max\,Supp\,}y^0, \ \ \mbox{or}\ \ \beta_0+\al_0\succeq{\rm max\,Supp\,}y^{\al_0}\end{equation}
and at least one of them with equality.
Note that \begin{equation}\label{al-be-gammm}
2(\al_0-\gamma) L_{\beta,1}=[L_{\beta-\al_0,0},L_{\al_0,1}]-[L_{\beta-\gamma,0},L_{\gamma,1}]
\mbox{ \ for \ }\beta,\gamma\in\G.\end{equation}First taking $\gamma=0$ and applying $\sigma$ to
it, by \eqref{sigma-l-a-0} and \eqref{max-be-----}
, we obtain
$y^\beta_\mu=0$ if $\mu\succ\beta_0+\beta$, i.e., \begin{equation}\label{suememe}
{\rm max\,Supp\,}y^\beta\preceq\beta_0+\beta\mbox{ \ for all \ }\beta\in\G,\end{equation} and
by computing the coefficient of $L_{\beta_0+\beta,j}$, we have (here and below, we simply denote $c^\beta_j=c^{\beta}_{\beta_0+\beta,j}$ for $\beta\in\G$),
\begin{equation}\label{cttttta}
c^{\beta}_j=\frac1{2\al_0}\Big((2\al_0+\beta_0-\beta)c^{\al_0}_{j}+(j-1)c^{\al_0}_{j-1}-
(\beta_0-\beta)c^{0}_{j}-(j-1)c^{0}_{j-1}\Big)=c_j+c'_j\beta,\end{equation}
where $c'_j=\frac1{2\al_0}(c^{0}_{j}-c^{\al_0}_{j}),$ $c_j=\frac1{2\al_0}\Big((2\al_0+\beta_0)c^{\al_0}_{j}+(j-1)c^{\al_0}_{j-1}-
\beta_0c^{0}_{j}-(j-1)c^{0}_{j-1}\Big)$ are some complex numbers not depending on $\beta$ but depending on $j$.
Applying $\sigma$ to \eqref{al-be-gammm} with arbitrary $\gamma$ and equating coefficients of $L_{\beta_0+\beta,j}$ gives
\begin{equation}\label{cttttta-gga}
2(\al_0-\gamma)c^{\beta}_j=(2\al_0+\beta_0-\beta)c^{\al_0}_{j}+(j-1)c^{\al_0}_{j-1}-
(2\gamma+\beta_0-\beta)c^{\gamma}_{j}-(j-1)c^{\gamma}_{j-1}.\end{equation}
Using \eqref{cttttta} to replace $c^{\beta}_j,c^{\gamma}_j$ and $c^{\gamma}_{j-1}$ in \eqref{cttttta-gga}, and comparing the coefficients of $\gamma^2$, we immediately obtain $c'_j=0$. Thus $c^\beta_j=c_j$ for all
$\beta\in\G$, and from \eqref{y-alpha}, we can write
\begin{equation}\label{l-al1-0}y^\beta_{\beta+\beta_0}=\sigma(L_{\beta,1})_{\beta+\beta_0}=\mbox{$\sum\limits_{j\ge1}$}c_jL_{\beta+\beta_0,j}.\end{equation}
Now apply $\sigma$ to
\begin{equation}\label{Applying-1to1}
L_{\al,2}=[L_{0,0},L_{\al,1}]-\al L_{\al,1},\end{equation}
 we see \begin{equation}\label{MSMAMA}
{\rm max\,Supp\,}\sigma(L_{\al,2})\preceq\al+\beta_0\mbox{ \  for \ }\al\in\G.\end{equation}
Applying $\sigma$ to \begin{equation}\label{Applying-2to2}
L_{\al,2}=\frac 1{\al-2\beta}[L_{\beta,1},L_{\al-\beta,1}]\mbox{ \ with \ }\al\ne2\beta,\end{equation} and computing the first term in its $(\al+2\beta_0)$-th homogeneous component, using \eqref{l-al1-0}, we see $\al+2\beta_0\in{\rm Supp\,}\sigma(L_{\al,2})$, a contradiction with \eqref{MSMAMA} and the fact that $\beta_0\succ0$. This proves Claim 2 (note that $c_1\ne0$ since $\sigma(W^1)=W^1\not\subset W^2$).\vskip5pt

Now from \eqref{Applying-1to1}, we obtain $\sigma(L_{\al,2})=\sum_{j\ge1}jc_jL_{\al,j+1}$.
Using this and \eqref{suememe}, \eqref{l-al1-0} (which is now equal to $y^{\beta}$), and applying $\sigma$ to \eqref{Applying-2to2}, we immediately obtain $c_1=1$, and $(k-1)c_{k-1}=\sum_{i+j=k}c_ic_j$ for $k\ge2$. But the latter relations force $c_k=c_2^{k-1}$ for all $k\geq 2$. So if $c_2\neq 0$, then $\sigma(L_{\alpha, 1})\notin W$, a contradiction.  Thus  $c_2=0$ and $\sigma(L_{\al,1})=L_{\al,1}$. Since $W$ is generated by $\{L_{\al,0},L_{\al,1}\,|\,\al\in\G\}$, we obtain $\sigma=1$.  This proves the theorem.\QED

\section{ Second cohomology group}
\setcounter{theo}{0}
Recall that a bilinear form $\psi:\L\times \L
\rightarrow \F$ is called a {\it 2-cocycle} on $\L$ if
the following conditions are satisfied:
\begin{eqnarray*}
&&\psi(x,y)=-\psi(y,x),
\ \ \ \psi(x,[y,z])+\psi(y,[z,x])+\psi(z,[x,y])=0,
\end{eqnarray*}
for $x,y,z\in \L$. Denote by
$C^{2}(\L,\ \F)$ the vector space of 2-cocycles on
 $\L$. For any linear function $f: \L \rightarrow
 \F$, one can define a 2-cocycles $\psi_{f}$ by
$ 
\psi_{f}(x,y)=f([x,y])$ for $x, y \in\L.$ 
Such a 2-cocycle is called a {\it trivial 2-cocycle} or a {\it 2-coboundary} on $\L$. Denote
by $B^{2}(\L,\ \F)$ the vector space of 2-coboundaries
on
 $\L$. The quotient space
 \begin{eqnarray*}
&&H^{2}(\L,\
\F)=C^{2}(\L,\ \F)/B^{2}(\L,\ \F)
\end{eqnarray*}
is called  the {\it 2-cohomology group} of $\L$.
There exists a one-to-one correspondence between the set of
equivalence classes of one-dimensional central extensions of
$\L$ by $\F$ and the 2-cohomology group of
$\L$.
\begin{theo}\label{theo4}We have $H^2(W,\C)=\C\ol{\phi}_0$, where $\ol{\phi}_0$ is the equivalence class of the 2-cocycle $\phi_0$, which is defined by
\begin{equation}\label{2-cococo}
\phi_0(L_{\al,i},L_{\beta,j})=\d_{\al+\beta,0}\d_{i+j,0}\frac{\al^3-\al}{12}.
\end{equation}
Thus $W$ has the unique universal central extension defined by \eqref{algebra-rela2}.
 \end{theo}
\noindent{\it Proof.~}~First note from \eqref{Aususus} that if $\G'=c\G$ for some $c\in\C^*$, then $W(\G)\cong W(\G')$. Thus without loss of generality, we can always suppose $1\in\G$.
Let $\psi\in C^2(W,\C)$, we define a linear function $f:W\to\C$ such that $f(L_{\al,i})$ is defined by induction on $i$ as follows
\begin{equation}\label{f====}
f(L_{\al,i})=\left\{\begin{array}{lll}
\frac1{2}\psi(L_{-1,0},L_{1,0})&\mbox{if \ }\al=i=0,\\[4pt]
\frac1\al\psi(L_{0,0},L_{\al,0})&\mbox{if \ }\al\ne0,\,i=0,\\[4pt]
\frac12\Big(\psi(L_{0,0},L_{0,1})+\psi(L_{-1,1},L_{1,0})\Big)&\mbox{if \ }\al=0,\,i=1,\\[4pt]
\frac1{2\al}\Big(\psi(L_{0,0},L_{\al,1})+\psi(L_{0,1},L_{\al,0})\Big)&\mbox{if \ }\al\ne0,\,i=1,\\[4pt]
\frac1{2}\Big(\psi(L_{0,0},L_{\al,1})-\psi(L_{0,1},L_{\al,0})\Big)&\mbox{if \ }i=2,\\[4pt]
\frac1{i-1}\Big(\psi(L_{0,0},L_{\al,i-1})-\al f(L_{\al,i-1})\Big)&\mbox{if \ }i\ge3.
\end{array}\right.
\end{equation}Set $\phi=\psi-\psi_f$. From the second and last cases of \eqref{f====}, we obtain
\begin{equation}
\label{phi===}
\phi(L_{0,0},L_{\al,i})=\psi(L_{0,0},L_{\al,i})-f([L_{0,0},L_{\al,i}])=0\mbox{ \ for \ }\al\in\G,\ 0\ne i\in\Z_+.
\end{equation}Similarly, from the second, fourth and fifth cases of  \eqref{f====}, we obtain
\begin{eqnarray}
\label{phi===++222}
\phi(L_{0,0},L_{\al,1})&\!\!\!=\!\!\!&\psi(L_{0,0},L_{\al,1})\!-\!f([L_{0,0},L_{\al,1}])\!=\!0,\\
\label{phi===++222+1}
\phi(L_{0,1},L_{\al,0})&\!\!\!=\!\!\!&\psi(L_{0,1},L_{\al,0})\!-\!f([L_{0,1},L_{\al,0}])\!=\!0
\mbox{ \ for \ }\al\in\G.
\end{eqnarray}
From the second and third cases of \eqref{f====},
\begin{equation}
\label{phi===++}
\phi(L_{1,0},L_{-1,1})=\psi(L_{1,0},L_{-1,1})-f([L_{1,0},L_{-1,1}])=0.
\end{equation}
Furthermore by replacing $\phi$ by $\phi-c\phi_0$ for some $c\in\C$, we can suppose
\begin{equation}
\label{phi===++1}
\phi(L_{\al,0},L_{\beta,0})=0\mbox{ \ for \ }\al,\beta\in\G.
\end{equation}
Then \eqref{phi===} gives\begin{equation}\label{Sec---}
0=\phi(L_{0,0},[L_{\al,0},L_{\beta,i-1}])=
(\al+\beta)\phi(L_{\al,0},L_{\beta,i-1})+(i-1)\phi(L_{\al,0},L_{\beta,i}).\end{equation}
Now we have
\begin{eqnarray}\label{A1111}
\!\!\!\!\!\!
({\al-2\beta})\phi(L_{\al,0},L_{\gamma,1})&\!\!\!=\!\!\!&
\phi([L_{\beta,0},L_{\al-\beta,0}],L_{\gamma,1})\nonumber\\
&\!\!\!=\!\!\!&(\beta-\gamma)\phi(L_{\al-\beta,0},L_{\beta+\gamma,1})
-\phi(L_{\al-\beta,0},L_{\beta+\gamma,2})\nonumber\\
&\!\!\!\!\!\!\!\!\!\!\!\!\!\!\!\!\!\!\!\!\!\!\!\!&
+
(\gamma+\beta-\al)\phi(L_{\beta,0},L_{\al-\beta+\gamma,1})
+\phi(L_{\beta,0},L_{\al-\beta+\gamma,2})\nonumber\\
&\!\!\!=\!\!\!&
(\al+\beta)\phi(L_{\al-\beta,0},L_{\beta+\gamma,1})
+(\beta-2\al)\phi(L_{\beta,0},L_{\al-\beta+\gamma,1}),
\end{eqnarray}where the last equality follows from \eqref{Sec---}.
Setting $\beta=-\al$ in \eqref{A1111} gives
\begin{eqnarray}\label{A1111+2}
\phi(L_{-\al,0},L_{\gamma+2\al,1})=-\phi(L_{\al,0},L_{\gamma,1}).\end{eqnarray}
Setting $\beta=-\gamma$ in \eqref{A1111}, by \eqref{phi===++222+1} and \eqref{A1111+2}, we obtain
\begin{eqnarray}\label{A1111+33}(\al+2\gamma)\phi(L_{\al,0},L_{\gamma,1})=(\gamma+2\al)\phi(L_{\gamma,0},L_{\al,1}).\end{eqnarray}
Exchanging $\al$ and $\gamma$ and changing $\beta$ to $-\beta$ in \eqref{A1111}, by \eqref{A1111+2} and \eqref{A1111+33}, we obtain
\begin{eqnarray}\label{A1111-------}
&\!\!\!\!\!\!&
(\gamma+2\beta)\frac{\al+2\gamma}{\gamma+2\al}\phi(L_{\al,0},L_{\gamma,1})\nonumber\\
&\!\!\!\!\!\!&=
(\gamma-\beta)\phi(L_{\gamma+\beta,0},L_{\al-\beta,1})+(\beta+2\gamma)\phi(L_{\beta,0},L_{\gamma+\al-\beta,1})\nonumber\\
&\!\!\!\!\!\!&=
(\gamma-\beta)\frac{\al+\beta+2\gamma}{2\al-\beta+\gamma}\phi(L_{\al-\beta,0},L_{\gamma+\beta,1})+(\beta+2\gamma)\phi(L_{\beta,0},L_{\gamma+\al-\beta,1})
,\end{eqnarray}
for all $\al,\beta,\gamma\in\G$ (note that if $\al=-2\gamma$ or  $\beta=2\al+\gamma$, we shall regard the equation as the one by first multiplying the equation by $\al+2\gamma$ or $2\al-\beta+\gamma$, then take $\al=-2\gamma$ or  $\beta=2\al+\gamma$). Taking $\beta=1$, we can solve $\phi(L_{\al,0},L_{\gamma,1})$ from \eqref{A1111} and \eqref{A1111-------} to obtain
\begin{equation}
\label{Soleee}
\phi(L_{\al,0},L_{\gamma,1})=\frac{\al\gamma(2\al+\gamma)}{(\al+\gamma+1)(\al+\gamma-1)}c_{\al+\gamma}\mbox{ \ if \ }\al\ne-\gamma\pm1,
\end{equation}
where $c_{\al}=\phi(L_{1,0},L_{\al-1,1})\in\C$. In particular, by \eqref{phi===++},
\begin{equation}
\label{Soleee-00}
\phi(L_{\al,0},L_{\gamma,1})=0\mbox{ \ \ if \ \ }\al+\gamma=0.
\end{equation}
Now \eqref{phi===} gives
\begin{equation}
\label{Soleee-0021}
0=\phi(L_{0,0},[L_{\al,1},L_{\beta,1}])=
(\al+\beta)\phi(L_{\al,1},L_{\beta,1})+\phi(L_{\al,2},L_{\beta,1})+\phi(L_{\al,1},L_{\beta,2}).
\end{equation}
Equation \eqref{Sec---} gives
\begin{eqnarray}\label{A1111-------++}
&\!\!\!\!\!\!\!\!\!\!\!\!\!\!\!\!\!\!\!\!\!\!\!\!&
-(\al+\beta)\beta\phi(L_{\al,0},L_{\beta,1})\nonumber\\
&\!\!\!\!\!\!\!\!\!\!\!\!\!\!\!\!\!\!\!\!\!\!\!\!&=\beta\phi(L_{\al,0},L_{\beta,2})
=\phi(L_{\al,0},[L_{0,1},L_{\beta,1}])\nonumber\\
&\!\!\!\!\!\!\!\!\!\!\!\!\!\!\!\!\!\!\!\!\!\!\!\!&
=
-\al\phi(L_{\al,1},L_{\beta,1})+\phi(L_{\al,2},L_{\beta,1})+(\beta-\al)\phi(L_{0,1},L_{\al+\beta,1})
+\phi(L_{0,1},L_{\al+\beta,2})
.\end{eqnarray}
Solving $\phi(L_{\al,2},L_{\beta,1})$ in this equation gives
\begin{eqnarray}\label{A1111---1++}
&&\phi(L_{\al,2},L_{\beta,1})\\
&=&
\al\phi(L_{\al,1},L_{\beta,1})+(\al-\beta)\phi(L_{0,1},L_{\al+\beta,1})
-\phi(L_{0,1},L_{\al+\beta,2})-(\al+\beta)\beta\phi(L_{\al,0},L_{\beta,1})\nonumber
.\end{eqnarray}
Thus we also have an expression for $\phi(L_{\al,1},L_{\beta,2})=-\phi(L_{\beta,2},L_{\al,1})$. Using them in
\eqref{Soleee-0021}, by \eqref{A1111+33}, we obtain
\begin{eqnarray}\label{A111-----++}
&\!\!\!\!\!\!&
\phi(L_{\al,1},L_{\beta,1})=\frac{\beta^2-\al^2}{2(2\al+\beta)}
\phi(L_{\al,0},L_{\beta,1})+\frac{\beta-\al}{\al+\beta}\phi(L_{0,1},L_{\al+\beta,1})\mbox{ \ if \ }\al\ne-\beta.
\end{eqnarray}
Using \eqref{Sec---}, we have
\begin{eqnarray}\label{A111---231--++}&&(\beta-\gamma)(\al+\beta+\gamma)\phi(L_{\al,0}, L_{\beta+\gamma,1}\nonumber)\\
&\!\!\!\!\!\!\!\!\!\!\!\!\!\!\!\!\!\!\!\!\!\!\!\!\!&=
(\gamma-\beta)\phi(L_{\al,0},L_{\beta+\gamma,2})\\
&\!\!\!\!\!\!\!\!\!\!\!\!\!\!\!\!\!\!\!\!\!\!\!\!\!&=
\phi(L_{\al,0},[L_{\beta,1},L_{\gamma,1}])
\nonumber\\
&\!\!\!\!\!\!\!\!\!\!\!\!\!\!\!\!\!\!\!\!\!\!\!\!\!&
=\!(\beta\!-\!\al)\phi(L_{\al+\beta,1},L_{\gamma,1})\!+\!\phi(L_{\al+\beta,2},L_{\gamma,1})
\!+\!(\gamma\!-\!\al)\phi(L_{\beta,1},L_{\al+\gamma,1})\!+\!\phi(L_{\beta,1},L_{\al+\gamma,2})\nonumber.
\end{eqnarray}
Using \eqref{Soleee},
\eqref{A1111---1++}
and \eqref{A111-----++} in \eqref{A111---231--++}, we can solve $c_\al=0$, thus
\eqref{Soleee}, \eqref{A111-----++} and \eqref{A1111---1++} give
\begin{eqnarray}\label{A111--21-++}
&\!\!\!\!\!\!&
\phi(L_{\al,0},L_{\beta,1})=0,\ \ \ \ \ \ \ \ \ \ \
\phi(L_{\al,1},L_{\beta,1})=\frac{\beta-\al}{\al+\beta}\tilde c_{\al+\beta}\mbox{ \ (if \ }\al\ne-\beta),\\
\label{A111--21-++2}
&\!\!\!\!\!\!&\phi(L_{\al,1},L_{\beta,2})=\frac{\al(\al-\beta)}{\al+\beta}\tilde c_{\al+\beta}+\bar c_{\al+\beta}\mbox{ \ (if \ }\al\ne-\beta),
\end{eqnarray}
where $\tilde c_{\al}=\phi(L_{0,1},L_{\al,1}),\,\bar c_\al=\phi(L_{0,1},L_{\al,2})\in\C$ (note that although we obtain the above under some conditions (for instance, there is a condition in \eqref{Soleee}), by \eqref{A1111}, all conditions can be removed as long as the expressions make senses).
Now using \eqref{A111--21-++2} in
$$(\gamma-\beta)\phi(L_{\al,1},L_{\beta+\gamma,2})=\phi(L_{\al,1},[L_{\beta,1},L_{\gamma,1}])=(\beta-\al)\phi(L_{\al+\beta,2},L_{\gamma,1})
+(\gamma-\al)\phi(L_{\beta,1},L_{\al+\gamma,2}),$$we obtain $\tilde c_\al=0$. Thus
\begin{eqnarray}\label{A111--21-1++}
&\!\!\!\!\!\!&
\phi(L_{\al,1},L_{\beta,1})=
0,\ \ \ \ \ \phi(L_{\al,1},L_{\beta,2})=\bar c_{\al+\beta}.
\end{eqnarray}
Using \eqref{Sec---} and \eqref{A111--21-++}, by induction on $i$, we obtain $\phi(L_{\al,0},L_{\beta,i})=0$. Thus
$$
0\!=\!\phi(L_{\al,0},[L_{\beta,1},L_{\gamma,2}])\!=\!(\beta\!-\!\al)\bar c_{\al+\beta+\gamma}\!+\!\phi(L_{\al+\beta,2},L_{\gamma,2})
\!+\!(\gamma\!-\!\al)\bar c_{\al+\beta+\gamma}\!+\!2\phi(L_{\beta,1},L_{\al+\gamma,3}).
$$
Setting $\al=-\beta$ and $\beta=0$ respectively, we can solve
\begin{eqnarray}\label{--21-1++}
&\!\!\!\!\!\!&
\phi(L_{\al,1},L_{\beta,3})=-\frac{(4\al+\beta)}{2}\bar c_{\al+\beta}-\frac12\hat c_{\al+\beta}
,\ \ \ \ \ \phi(L_{\al,2},L_{\beta,2})=3\al\bar c_{\al+\beta}+\hat c_{\al+\beta},
\end{eqnarray}
where $\hat c_\al=\phi(L_{0,2},L_{\al,2}).$ Since $ \phi(L_{\al,2},L_{\beta,2})=-\phi(L_{\beta,2},L_{\al,2})$, we obtain $\hat c_\al=-\frac32\al\bar c_\al$.
Using \eqref{--21-1++} in
$$\begin{array}{lll}
(\beta-\al)\phi(L_{0,1},L_{\al+\beta,3})+\phi(L_{0,1},L_{\al+\beta,4})&\!\!\!=\!\!\!&\phi(L_{0,1},[L_{\al,1},L_{\beta,2}])
\\[4pt]
&\!\!\!=\!\!\!&\al\phi(L_{\al,2},L_{\beta,2})+\beta\phi(L_{\al,1},L_{\beta,3})+\phi(L_{\al,1},L_{\beta,4}),\end{array}$$
we solve \begin{equation}\label{y14==}
\phi(L_{\al,1},L_{\beta,4})
=\frac{\al}4(11\beta-7\al)\bar c_{\al+\beta}+c'_{\al+\beta},\end{equation}
 where
$c'_\al=\phi(L_{0,1},L_{\al,4})\in\C.$
Using this and \eqref{--21-1++} in
$$
0=\phi(L_{0,0},[L_{\al,1},L_{\beta,3}])=(\al+\beta)\phi(L_{\al,1},L_{\beta,3})
+\phi(L_{\al,2},L_{\beta,3})+3\phi(L_{\al,1},L_{\beta,4}),
$$
we obtain
\begin{equation}\label{y23==}\phi(L_{\al,2},L_{\beta,3})=
\frac14(26\al^2-29\al\beta-\beta^2)
\bar c_{\al+\beta}-3c'_{\al+\beta}.\end{equation}
Using this and \eqref{--21-1++} in
$$0=\phi(L_{0,0},[L_{\al,2},L_{\beta,2}])=
(\beta+\al)\phi(L_{\al,2},L_{\beta,2})+2\phi(L_{\al,3},L_{\beta,2})+
2\phi(L_{\al,2},L_{\beta,3})
,$$
we immediately obtain $\bar c_\al=0$. This proves $$\phi(L_{\al,i},L_{\beta,j})=0\mbox{ \ \ for \ \ }i+j\le4.$$
Now inductively assume for $k\ge4$, we have proved $\phi(L_{\al,i},L_{\beta,j})=0$ for all $i+j\le k$. Then
\begin{eqnarray}\label{i1===}
0=\phi(L_{\al,0},[L_{\beta,i},L_{\gamma,k-i}])=i\phi(L_{\al+\beta,i+1},L_{\gamma,k-i})
+(k-i)\phi(L_{\beta,i},L_{\al+\gamma,k-i+1}).
\end{eqnarray}
Setting $\beta=0$ and $\gamma=0$ in \eqref{i1===} respectively, we can solve
\begin{equation}\label{FIEL}
\phi(L_{\al,i},L_{\beta,k+1-i})=(-1)^{i-1}\binom{k-1}{i-1}c^k_{\al+\beta}\mbox{ \ for some \ }c^k_\al\in\C.\end{equation}
If $k=2\ell-1$ is odd, then since $\phi(L_{\al,\ell}, L_{\beta,\ell})=-\phi(L_{\beta,\ell}, L_{\al,\ell})$, we obtain $c^k_\al=0$ by \eqref{FIEL}, and thus $\phi(L_{\al,i},L_{\beta,k+1-i})=0$. Now assume $k=2\ell$ is even.
 Using \eqref{FIEL}, we have
\begin{eqnarray}\label{last2}
0&\!\!\!=\!\!\!&\phi(L_{\al,0},[L_{\beta,\ell},L_{\gamma,\ell+1}])
\nonumber\\
&\!\!\!=\!\!\!&(\beta-\al)\phi(L_{\al+\beta,\ell},L_{\gamma,\ell+1})
+\ell\phi(L_{\al+\beta,\ell+1},L_{\gamma,\ell+1})
\nonumber\\&\!\!\!\!\!\!&+(\gamma-\al)\phi(L_{\beta,\ell},L_{\al+\gamma,\ell+1})
+(\ell+1)\phi(L_{\beta,\ell},L_{\al+\gamma,\ell+2})
\nonumber\\
&\!\!\!=\!\!\!&
(\beta+\gamma-2\al)\phi(L_{0,\ell},L_{\al+\beta+\gamma,\ell+1})
\nonumber\\&\!\!\!\!\!\!&+\ell\phi(L_{\al+\beta,\ell+1},L_{\gamma,\ell+1})
+(\ell+1)\phi(L_{\beta,\ell},L_{\al+\gamma,\ell+2}).
\end{eqnarray}
Setting $\beta=0$, we obtain
\begin{eqnarray}\label{last2++}
0&\!\!\!=\!\!\!&
(\gamma-2\al)\phi(L_{0,\ell},L_{\al+\gamma,\ell+1})
+\ell\phi(L_{\al,\ell+1},L_{\gamma,\ell+1})+(\ell+1)\phi(L_{0,\ell},L_{\al+\gamma,\ell+2}).
\end{eqnarray}
Exchanging $\al$ and $\gamma$, and summing the result with \eqref{last2++},
by the skew symmetry of \linebreak $\phi(L_{\al,\ell+1},L_{\gamma,\ell+1})$, we  \vspace*{-7pt}obtain
$$\phi(L_{0,\ell},L_{\al,\ell+2})=\frac{\al}{2(\ell+1)}\phi(L_{0,\ell},L_{\al,\ell+1})
.$$
Using this in \eqref{last2++} gives $\phi(L_{\al,\ell+1},L_{\gamma,\ell+1})=\frac{3(\gamma-\al)}{2\ell}\phi(L_{0,\ell},L_{\al+\gamma,\ell+1})$.
Using this in \eqref{last2}, we obtain $$0=\frac{5\gamma-\beta-7\al}{2}
\phi(L_{0,\ell},L_{\al+\beta+\gamma,\ell+1})
+(\ell+1)\phi(L_{\beta,\ell},L_{\al+\gamma,\ell+2}).$$
Exchanging $\al$ and $\gamma$ shows $\phi(L_{0,\ell},L_{\al+\beta+\gamma,\ell+1})=0$.
This together with \eqref{FIEL} proves $c^k_\al=0$, and thus
$\phi(L_{\al,i},L_{\beta,k+1-i})=0$. The theorem is proved by induction.\QED



\end{CJK*}
\end{document}